\documentclass[a4paper]{article}

\usepackage{graphicx}
\usepackage{latexsym}
\usepackage{amsmath}
\usepackage{amssymb}
\usepackage{caption} 

\newenvironment{itemise}{\begin{itemize}}{\end{itemize}}
\newenvironment{centre}{\begin{center}}{\end{center}}

\newtheorem{definition}{Definition}
\newtheorem{proposition}{Proposition}

\newcommand{\paren}[1]{\ensuremath{\left(#1\right)}}

\newcommand{\pwd}{\ensuremath{\mathrm{PWD}}}
\newcommand{\dwp}{\ensuremath{\mathrm{DWP}}}
\newcommand{\od}{\ensuremath{\mathrm{OD}}}
\newcommand{\ga}{\ensuremath{A_0}}
\newcommand{\gp}{\ensuremath{P_0}}
\newcommand{\ua}{\ensuremath{A_U}}
\newcommand{\up}{\ensuremath{P_U}}
\newcommand{\pd}{\ensuremath{D}}
\newcommand{\uai}{\ensuremath{A_U^{(1)}}}
\newcommand{\upi}{\ensuremath{P_U^{(1)}}}
\newcommand{\pdi}{\ensuremath{D^{(1)}}}
\newcommand{\pwdi}{\ensuremath{\mathrm{PWD}^{(1)}}}
\newcommand{\uaii}{\ensuremath{A_U^{(2)}}}
\newcommand{\upii}{\ensuremath{P_U^{(2)}}}
\newcommand{\pdii}{\ensuremath{D^{(2)}}}
\newcommand{\pwdii}{\ensuremath{\mathrm{PWD}^{(2)}}}
\newcommand{\subdiv}{\ensuremath{\mathfrak{S}}}

\begin{document}

\title{Population-Weighted Density, Density-Weighted Population, Granularity, Paradoxes: a Recapitulation}
\author{Anthony B.\ Morton\\Honorary Fellow, University of Melbourne, Australia}
\date{January 2015}
\maketitle

\section{Introduction}

Quantifying the population density of an urban area is a fraught issue, related closely to that of measuring trends in popular but imprecise urban development concepts such as `smart growth' and `sprawl'.
Measures of density are often defined differently from place to place or inconsistently applied. Arguments abound meanwhile over just how much of the land surrounding a city should and should not be classified as urban\footnote{`Urban' in this note is a synonym for the built-up or urbanised land making up a metropolitan area, as distinct from rural or interurban land.
Especially in North American contexts, `urban' can be given a more narrow connotation in opposition to `suburban'.
Thia latter concept will generally be termed `intense urban' or `inner-urban' in this note.}, and how (or indeed whether) one may draw a distinction between `inner-urban' and `suburban' parts of a metropolitan area.

The most popular and straightforward method for calculating urban population density is to draw a notional boundary around the urban area to be measured, and then simply divide the population within the boundary by the total area of land enclosed.
This gives the quantity defined below as \emph{overall density}, also known as \emph{average density} or \emph{gross density}.
It is however a contentious measure, owing to the question of where to place the boundary.
The pitfalls and gross misconceptions that can arise from naive calculations using arbitrary administrative boundaries, which often bear little relation to the actual urbanised area, were noted as early as 1946 by the Viennese--Australian town planner and architect Dr Ernest Fooks \cite{fooks1946} and more recently by Mees \cite{mees2000,mees2010}.

A number of more refined approaches have been developed in response to this problem.
One of the more well-established methods, foreshadowed by Fooks but developed further by Linge \cite{linge1965} and others, retains the basic definition but takes a more rigorous approach to the boundary.
For this, the region is divided into small \emph{parcels} (given various names such as `census tracts', `collection districts' or `statistical areas') and the unique boundary drawn that separates parcels above a fixed `cutoff density' on one side from parcels below it on the other.
This cutoff density is chosen sufficiently low that most regions fitting the everyday notion of a `built-up area' will comfortably exceed it: 2 persons per hectare\footnote{Land area is measured in hectares or acres depending on local convention.  The distinction is not important for the purpose of this discussion, and hectares are used for concreteness.  The conversions 1 hectare $\approx$ 2.5 acres and 1000 hectares $\approx$ 4 square miles may be used as first approximations.  100 hectares is 1 square kilometre.} is the value used by the Australian Bureau of Statistics, following \cite{linge1965}.
As many European, North American and Australian cities have fairly well-defined contiguous urban areas, this has provided a relatively useful measure.

Other methods of quantifying urban density also anticipated by Fooks include calculating overall density within concentric rings at set distances from a defined urban centre, overlaying a grid and indicating the overall density for each grid element in a diagram, or computing a Gini coefficient that quantifies the `unevenness' of the population distribution within the urban area.
Eidlin \cite{eidlin2010} compares a number of these measures while seeking to quantify the popular concept of \emph{urban sprawl}.
Eidlin is particularly captivated by the `Los Angeles paradox'---the fact that although Los Angeles is widely recognised as a highly car-dependent and spread-out city, this is not well reflected in actual measures of urban density.
Indeed on the US Census Bureau's overall density measure (based on a refined approach similar to Linge's) Los Angeles is the most dense urban region in the USA.
This stems from the fact that while LA lacks the dense urban cores of cities such as New York or San Francisco, its suburbs have been developed at densities greater than typical of other US suburbs, and with relatively little variation across the urban area \cite{eidlin2010}.
This medium-density, suburb-dominated pattern is also typical of Australian cities \cite{mees2000}.

The measure known as \emph{population-weighted density} (PWD, also called \emph{perceived density}) is one that has recently emerged as a challenge to older definitions such as overall density \cite{bradford2008}.
The US Census Bureau has adopted PWD as a density measure since its 2010 Census \cite{uscb2012}.
As one of its original proponents Chris Bradford points out, PWD does not completely resolve the Los Angeles paradox, with Los Angeles still ranking in PWD terms ahead of every other US city except New York and San Francisco \cite{bradford2008,eidlin2010}.
Nonetheless, it is convincingly argued that PWD provides a superior measure to overall density because it ``gives equal weighting to each resident, rather than to each hectare of land''.
It is also held to be less sensitive to the detail of how one delineates an urban area, as the definition of PWD automatically discounts sparsely populated non-urban land.

In this note some less intuitive aspects of PWD are explored, so that the consequences of adopting PWD as a density measure are better understood.
In particular, it is recalled (as the definition in \cite{uscb2012} implies) that PWD bears a close relationship to the apparently inverted concept of a \emph{density-weighted population} (DWP).
It will also be seen that one cannot entirely dispense with the need to define urban boundaries, or to work preferentially with the smallest parcels of land for which one has data.

\section{The Key Formula}

Population-weighted density is defined as follows.
Start with a conveniently defined region of area \ga\ (the \emph{gross area}) and population \gp, which contains the entire urban area in question.
Now divide this up into smaller parcels of land: let the area of the $k$th parcel be $A_k$ and its population $P_k$.
The average density of the $k$th parcel is $P_k / A_k$.
The PWD is a weighted sum, where each parcel density is weighted by the parcel's share $P_k / \gp$ of the total population.
So, if there are $N$ parcels in total:
\begin{equation}
\pwd = \sum_{k=1}^N \frac{P_k}{\gp} \cdot \frac{P_k}{A_k}.
\label{eq:pwd}
\end{equation}
One may compare formula (\ref{eq:pwd}) with that for overall density, which is
\begin{equation}
\od = \frac{\gp}{\ga} = \frac{\sum_{k=1}^N P_k}{\sum_{k=1}^N A_k}
   = \sum_{k=1}^N \frac{A_k}{\ga} \cdot \frac{P_k}{A_k}.
\label{eq:od}
\end{equation}
The rightmost of these equivalent expressions shows why in comparison with PWD, OD is sometimes referred to as `area-weighted density'---thus motivating the concept of PWD as a measure based on `people rather than hectares'.

Now the alternative concept of density-weighted population (DWP) is defined.
Again, it is a weighted sum based on division into parcels with population $P_k$ and area $A_k$.
Now, however, one takes the population $P_k$ of each parcel and weights it according to the parcel's relative density: the ratio of its actual density $P_k / A_k$ to the overall density for the entire region given by (\ref{eq:od}).
Summing over all parcels, this gives
\begin{equation}
\dwp = \sum_{k=1}^N P_k \cdot \frac{P_k / A_k}{\od}.
\label{eq:dwp}
\end{equation}
Thus, if a given parcel of land is twice as dense in population as the region as a whole, each resident of that parcel counts as 2 people in the DWP measure.

Although PWD and DWP are defined rather differently and appear to measure different concepts, they are in fact almost the same thing.
The key formula linking the two is
\begin{equation}
\pwd = \frac{\dwp}{\ga}.
\label{eq:pwddwp}
\end{equation}
This identity follows immediately from the definitions (\ref{eq:pwd}), (\ref{eq:od}) and (\ref{eq:dwp}).
The explicit calculation is
\begin{equation}
\frac{\dwp}{\ga} = \frac{1}{\ga} \sum_{k=1}^N P_k \cdot \frac{P_k / A_k}{\gp / \ga}
   = \sum_{k=1}^N \frac{P_k}{\gp} \cdot \frac{P_k}{A_k} = \pwd.
\label{eq:pwddwpc}
\end{equation}
One way to view the formula (\ref{eq:pwddwp}) is as follows.
As with all measures of population density, the units of PWD are persons per hectare.
The `hectare' dimension is supplied by the gross area \ga: the `person' dimension, however, is supplied by a density-weighted population rather than \gp\ itself.

Accordingly, PWD is in a sense the average density when residents themselves are weighted according to the relative density of their neighbourhoods.
It is indeed a measure based on people, but it is less clear that people receive more equal treatment in this measure, as is often claimed.

\section{PWD and Subdivision of Parcels}

Viewing population-weighted density PWD via a weighted population DWP helps in developing an intuitive picture of how PWD numbers change when large parcels of land are subdivided into smaller parcels.
Unlike measures such as overall density, which depend only on how the boundary of the urban area is defined, the PWD is sensitive to the size and distribution of the small parcels.

As will be shown formally in the next section, PWD never falls when land is divided into smaller parcels, and almost always increases.
But there is a broad spectrum of behaviour depending on how much the density varies within a region.

At one end of the spectrum, if density within the region is absolutely uniform, then the PWD is equal to the overall density no matter what subdivision is used; the densities being weighted are all identical and so the weighting has no effect.
It is similarliy intuitive that in this case DWP is always equal to the regional population \gp: this is because the population of each parcel is weighted each time by the same relative density 1.

At the opposite end of the spectrum is where the population \gp\ is concentrated within a single parcel of area $A_1$---or more generally, within $N$ parcels each of the same density and with combined area $A_1$.
The remaining area $(\ga - A_1)$ is unoccupied.
In this case, the PWD is equal to the density $\gp / A_1$ of the concentrated population.
Notice that when this is found by applying the DWP formula, it results from recognising that the relative density of a populated land parcel is $\ga / A_1$, the ratio of the gross area to the occupied area.

Of more interest are the intermediate cases where the small parcels are of varying density, or where the density actually varies significantly within a single parcel.
The use of overly large parcels with varying internal density actually has paradoxical effects as an urban area evolves in time, as will be shown in a later section.

The alternative route to PWD calculations via the DWP formula can aid intuition in many cases.
For example, suppose a region with variable population density is subdivided into two parcels, one with twice the density of the wider region, and the other with half the density.
If the parcels each have the same population $\gp / 2$, then the DWP calculation is
\begin{equation}
\dwp = \paren{2} \paren{\frac{\gp}{2}} + \paren{\frac{1}{2}} \paren{\frac{\gp}{2}} = \frac{5}{4} \gp
\label{eq:subex1a}
\end{equation}
and the resulting PWD is
\begin{equation}
\pwd = \frac{5}{4} \frac{\gp}{\ga} = \frac{5}{4} \od
\label{eq:subex1b}
\end{equation}
where \od\ is the overall density of the original region.
In other words, this choice of subdivision has effectively increased both the DWP and PWD figures by 25 per cent.

Or suppose the region is divided into three parcels of equal population $\gp / 3$, where two parcels have density twice that of the overall region, and the third, half the density.
Again the DWP formula is straightforward:
\begin{equation}
\dwp = \paren{2} \paren{\frac{\gp}{3}} + \paren{2} \paren{\frac{\gp}{3}}
   + \paren{\frac{1}{2}} \paren{\frac{\gp}{3}} = \frac{3}{2} \gp
\label{eq:subex2a}
\end{equation}
and
\begin{equation}
\pwd = \frac{3}{2} \frac{\gp}{\ga} = \frac{3}{2} \od.
\label{eq:subex2b}
\end{equation}
Here, the effective increase in DWP or PWD is 50 per cent.

As a general rule one sees that when a larger area \ga\ is subdivided, the ratio $\pwd / \od$ of weighted density to overall density (what Eidlin \cite{eidlin2010} calls the \emph{density gradient index}) is identical to the ratio $\dwp / \gp$ of density-weighted to ordinary population.

\section{Subdivision and Harmonic Means}

In this section some more formal mathematical results are presented which underlie the informal discussion of the previous section.

Population-weighted density has a close connection with \emph{harmonic means}.
Given a set of values $x_1, x_2, \ldots, x_N$, the harmonic mean is the number $x_H$ such that
\begin{equation}
\frac{1}{x_H} = \frac{1}{N} \paren{\frac{1}{x_1} + \frac{1}{x_2} + \ldots + \frac{1}{x_N}}.
\label{eq:hmean}
\end{equation}
Note that if all the $x_k$ are identical, then $x_H$ also takes the same value.
More generally, given any set of weights $w_1, w_2, \ldots, w_N$ one may define a \emph{weighted harmonic mean} $x_H^{(w)}$ by the formula
\begin{equation}
\frac{1}{x_H^{(w)}} = \frac{1}{w_1 + w_2 + \ldots + w_N}
   \paren{\frac{w_1}{x_1} + \frac{w_2}{x_2} + \ldots + \frac{w_N}{x_N}}.
\label{eq:whmean}
\end{equation}
Of course, if all the $w_k$ are equal then $x_H^{(w)}$ reduces to the ordinary harmonic mean.

There are well-known inequalities relating both the ordinary and weighted harmonic means to the corresponding arithmetic means.
The \emph{weighted arithmetic mean} of the $x_k$ with weights $w_k$ is
\begin{equation}
\bar{x}^{(w)} = \frac{1}{w_1 + w_2 + \ldots + w_N}
   \paren{w_1 x_1 + w_2 x_2 + \ldots + w_N x_N}.
\label{eq:wamean}
\end{equation}
If all the weights are equal then $\bar{x}^{(w)}$ is the ordinary arithmetic mean, denoted $\bar{x}$.
The \emph{weighted power means inequality} asserts (among other things) that
\begin{equation}
x_H^{(w)} \leq \bar{x}^{(w)} \qquad \text{and specifically} \qquad x_H \leq \bar{x}
\label{eq:ahmeanie}
\end{equation}
with equality if and only if all the $x_k$ are equal.

In the case of the ordinary means, a lower bound is known for the difference $\bar{x} - x_H$, in the case where the $x_k$ are all positive.
If the variance of the value set $\{x_1, x_2, \ldots, x_N\}$ is $\sigma^2$, and the largest of the $x_k$ is no greater than $M$, then one has
\begin{equation}
\bar{x} - x_H \geq \frac{\sigma^2}{2 M}.
\label{eq:ahbound}
\end{equation}

Now, suppose a region is subdivided into $N$ parcels, each containing an \emph{equal} population $\gp / N$.
Then each parcel contains an equal share $1 / N$ of the population, and the PWD is
\begin{equation}
\pwd = \frac{1}{N} \paren{\frac{\gp / N}{A_1} + \frac{\gp / N}{A_2} + \ldots + \frac{\gp / N}{A_N}}
   = \frac{\gp}{N^2} \paren{\frac{1}{A_1} + \frac{1}{A_2} + \ldots + \frac{1}{A_N}}
   = \frac{\gp}{N A_H}
\label{eq:pwdeq}
\end{equation}
where $A_H$ is the harmonic mean of the parcel areas $A_k$.
However, by the power-means inequality (\ref{eq:ahmeanie}) one has
\begin{equation}
A_H \leq \bar{A} = \frac{A_1 + A_2 + \ldots + A_N}{N} = \frac{\ga}{N}
\label{eq:ahie}
\end{equation}
and hence
\begin{equation}
\pwd = \frac{\gp}{N A_H} \geq \frac{\gp}{N (\ga / N)} = \frac{\gp}{\ga} = \od
\label{eq:pwdpdie}
\end{equation}
with equality only if all parcels are the same size---equivalently, the same density.
One can in fact bound the difference between PWD and OD (or rather their reciprocals) using the inequality (\ref{eq:ahbound}), which by virtue of (\ref{eq:pwdeq}) and (\ref{eq:ahie}) becomes
\begin{equation}
\frac{\ga}{N} - \frac{\gp}{N \cdot \pwd} \geq \frac{\sigma^2}{2 M}
\label{eq:pwdbound1}
\end{equation}
or
\begin{equation}
\frac{1}{\od} - \frac{1}{\pwd} \geq \frac{N \sigma^2}{2 M \gp}
\label{eq:pwdbound2}
\end{equation}
where $\sigma$ is the standard deviation of the parcel areas, and $M$ is the area of the largest parcel.

The result $\pwd \geq \od$ extends to the more general case where the subdivision is into parcels with unequal populations $P_k$.
However, the (weighted) harmonic mean should now be formed not from the areas $A_k$ but from the inverse densities, $\Delta_k = A_k / P_k$.
Express the PWD calculation as
\begin{equation}
\pwd = \frac{1}{\gp} \paren{P_1 \frac{P_1}{A_1} + P_2 \frac{P_2}{A_2} + \ldots + P_N \frac{P_N}{A_N}}
   = \frac{1}{\Delta_H^{(P)}}
\label{eq:pwdneq}
\end{equation}
where $\Delta_H^{(P)}$ is a weighted harmonic mean of the inverse densities, with weights equal to the populations $P_k$.
The associated weighted arithmetic mean is
\begin{equation}
\bar{\Delta}^{(P)}
   = \frac{1}{\gp} \paren{P_1 \frac{A_1}{P_1} + P_2 \frac{A_2}{P_2} + \ldots + P_N \frac{A_N}{P_N}}
   = \frac{\ga}{\gp} = \frac{1}{\od}
\label{eq:bardeltap}
\end{equation}
and so the result corresponding to (\ref{eq:pwdpdie}) with unequal populations is
\begin{equation}
\pwd = \frac{1}{\Delta_H^{(P)}} \geq \frac{1}{\bar{\Delta}^{(P)}} = \od
\label{eq:pwdpdiegen}
\end{equation}
with equality only when the quantities $\Delta_k$ are equal---that is, when all the parcels have the same density, just as in the case of equal populations.

There remains one technical gap to deal with: while one does not lose anything of consequence by requiring that $A_k > 0$, it may sometimes come about that $P_k = 0$ for one or more parcels, in which case $\Delta_k$ is undefined.
In this case, let $A_Z \geq 0$ denote the total area of parcels with zero population, and apply the above argument to the subregion of area $\ga - A_Z$ comprising all those parcels with $P_k > 0$.
The PWD of this subregion is equal to the PWD of the original region, since in the latter all the parcels with zero population have zero weighting.
But by (\ref{eq:pwdpdiegen}) this PWD is also greater than or equal to the overall density of the subregion, which is $\gp / (\ga - A_Z) \geq \gp / \ga = \od$.
So it remains true that $\pwd \geq \od$, even when some of the parcels may have zero population.
(Indeed if there actually are parcels with zero population then $A_Z > 0$ and the inequality is strict: $\pwd > \od$.)

One can collect all the above into the following formal result:
\begin{proposition}
\label{prop:pwdpd}
Let a single region of population \gp\ and gross area \ga\ be subdivided into any number $N$ of smaller parcels, each of population $P_k$ and area $A_k$ with $A_k > 0$.
Then the PWD of the region based on the $N$ parcels is always greater than or equal to the overall density $\od = \gp / \ga$.
It equals \od\ if and only if the population densities $P_k / A_k$ of the parcels are all equal to one another (and hence to \od).
In the case where all parcels have the same population, one has the lower bound
\begin{equation}
\frac{1}{\od} - \frac{1}{\pwd} \geq \frac{N \sigma_A^2}{2 A_{\max} \gp}
   = \frac{\sigma_A}{A_{\max}} \cdot \frac{\sigma_A / 2}{\gp / N}
\label{eq:pwdbound}
\end{equation}
where $\sigma_A$ is the standard deviation of the parcel areas and $A_{\max}$ the largest parcel area.
\end{proposition}

Now suppose that a larger region has already been subdivided into parcels, and the $k$th parcel (say) is further subdivided into smaller parcels.
Proposition \ref{prop:pwdpd} applies to the subdivision of parcel $k$, just as it does to the original region, so one may assert that $\pwd_k$---the PWD of parcel $k$ alone based on this second subdivision---is greater than or equal to the overall density $P_k / A_k$ of parcel $k$, with equality only if all the smaller parcels also have density $P_k / A_k$.

Now let the PWD of the original region be calculated in two ways: the first treating parcel $k$ as one parcel, using the original subdivision; the second using the further subdivision of parcel $k$ into smaller parcels, with the other parcels unchanged.
These two calculations will differ precisely as follows: a single term $(P_k / \gp) (P_k / A_k)$ in the first calculation is replaced with $(P_k / \gp) \pwd_k$ in the second calculation.
(In the second case, the factor $P_k / \gp$ is required in order to change the normalisation of the terms in $\pwd_k$.)
But since these are purely additive terms in the larger PWD calculation, and $\pwd_k \geq P_k / A_k$, it follows that the PWD calculated the second way---using the finer division for parcel $k$---is always greater than or equal to the PWD calculated the first way (again, with equality only if all the parcel densities are equal).

The same argument may now be repeated for another parcel other than parcel $k$, and again as many times as desired, until all the original parcels have been subdivided.
At every stage of the argument, it remains true that the PWD based on the finer subdivision is no less than the PWD for the coarser one, and equal only if all the `before and after' parcel densities are the same.

To state the final result formally requires some more precise definitions:
\begin{definition}
\label{def:subdiv}
A \emph{subdivision} \subdiv\ of a region $R$ (viewed as a compact two-dimensional point set) is a covering of $R$ by some number $N$ of (closed) subsets $S_k$, $1 \leq k \leq N$, each of nonzero area, such that $S_k \cap S_j$ has zero area whenever $k \neq j$, and $\cup_k S_k = R$.
If $\subdiv_1$ and $\subdiv_2$ are two subdivisions of the same region $R$, then $\subdiv_2$ is a \emph{proper subdivision} of $\subdiv_1$ if every subset $S_k \in \subdiv_1$ is equal to the union of some collection of subsets $S_i', S_j', \ldots \in \subdiv_2$.
\end{definition}
The references to `compact' and `closed' sets are mathematical technicalities for the sake of precision: they amount to requiring that every set has a well-defined boundary.
It is worth noting that any region $R$ has a `trivial' subdivision $\subdiv = \{R\}$, in which $R$ acts as a subdivision of itself.
As another `trival' fact, one may observe that any subdivision is a proper subdivision of itself.
Last but not least, observe that when $\subdiv_2$ is a proper subdivision of $\subdiv_1$, every set $S_k \in \subdiv_2$ is contained in a unique set $S_k' \in \subdiv_1$.
\begin{proposition}
\label{prop:pwdsubdiv}
Let $\subdiv_1$ and $\subdiv_2$ be two subdivisions of a region $R$, such that $\subdiv_2$ is a proper subdivision of $\subdiv_1$.
Let $\pwd_1$ be the PWD calculated for $R$ using the subset populations $P_k$ and areas $A_k$ of the subdivision $\subdiv_1$, and $\pwd_2$ be calculated likewise using subdivision $\subdiv_2$.
Then one has
\begin{equation}
\pwd_2 \geq \pwd_1
\label{eq:pwdsubdiv}
\end{equation}
with equality only if every set $S_k \in \subdiv_2$ has the same population density as the set $S_k' \in \subdiv_1$ that contains $S_k$.
\end{proposition}
Note that if for $\subdiv_1$ one takes the trivial subdivision $\{R\}$, Proposition \ref{prop:pwdsubdiv} reduces to Proposition \ref{prop:pwdpd} (excluding the lower bound result, which relies on equal populations $P_k$).

\section{The Perils of Large Parcels I: When Density Increase Leads to Falling PWD}

An evident advantage of the PWD (\ref{eq:pwd}) over overall density (\ref{eq:od}) is that including non-urban land within the area \ga\ does not greatly affect the PWD quantity, since any parcels comprising \emph{wholly} non-urban land will have small population and therefore receive negligible weighting in the calculation.
Whereas in (\ref{eq:od}), it is important that the region itself be chosen to exclude non-urban land, otherwise the density number obtained will severely underestimate the true urban density.

The automatic discounting of non-urban land by (\ref{eq:pwd}), however, breaks down when the parcels are sufficiently large to contain significant amounts of both urban \emph{and} non-urban land.
In this case, PWD calculations can actually have paradoxical results: in particular, an urban area that expands into its non-urban hinterland can appear to have a declining PWD based on a larger enclosing region, \emph{even when the density of the urbanised area remains constant or increases modestly}.

This can be demonstrated using an idealised example with just two parcels: an `inner' parcel entirely within the urban area, and an `outer' parcel comprising both urban and non-urban land.
The example is illustrated schematically in Figure \ref{fig:schematic} and constructed as follows:
\begin{figure}
\begin{centre}
\includegraphics[width=15cm]{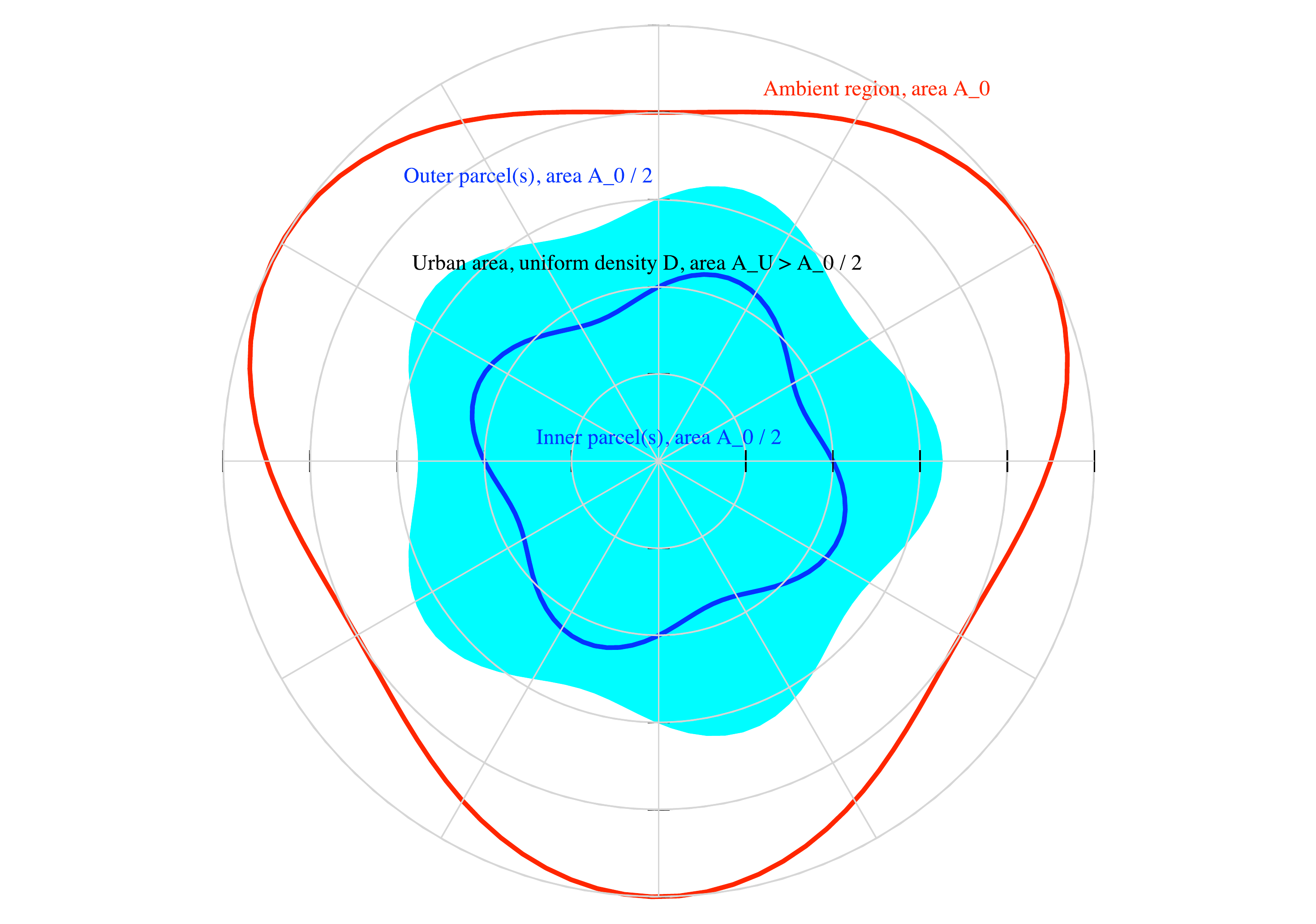}
\end{centre}
\caption{Schematic representation of urban growth example.}
\label{fig:schematic}
\end{figure}
\begin{itemise}
\item
Both parcels have the same area $\ga / 2$.
\item
The actual urbanised area of the city is $\ua < \ga$ and increases with time.
\item
The inner parcel is entirely included within the area \ua (and hence $\ua > \ga / 2$ at all times).
\item
The population density $\pd = \up / \ua$ within the urban area also varies with time, but at any given point in time is the same in all parts of the urban area \ua.
\item
There is no population outside the urban area, so $\up = \gp$ at all times, and it may be said that every resident within the entire region experiences the same population density \pd.
\end{itemise}

Now consider two snapshots in time: an earlier snapshot with urban area, population and density equal to \uai, \upi\ and \pdi\ respectively, and a later snapshot where these are given by \uaii, \upii\ and \pdii\ respectively, with $\uaii > \uai$.

At the earlier time, the inner parcel of land contains a proportion $(\ga / 2) / \uai$ of the actual urban area, and since the urban density is uniform, the population of each parcel is given as
\begin{equation}
P_1 = \frac{\ga}{2 \uai} \upi, \qquad P_2 = \paren{1 - \frac{\ga}{2 \uai}} \upi.
\label{eq:p12}
\end{equation}
The parcel areas are $A_1 = A_2 = \ga / 2$, and so the densities of each parcel are
\begin{equation}
\frac{P_1}{A_1} = \frac{\upi}{\uai} \paren{= \pdi}, \qquad
\frac{P_2}{A_2} = \frac{2 \upi}{\ga} - \frac{\upi}{\uai}.
\label{eq:d12}
\end{equation}
Putting these into the formula for PWD, one has
\begin{align}
\pwdi &= \frac{\ga}{2 \uai} \cdot \frac{\upi}{\uai}
      + \paren{1 - \frac{\ga}{2 \uai}} \paren{\frac{2 \upi}{\ga} - \frac{\upi}{\uai}} \nonumber \\
   &= \frac{\ga}{\uai} \cdot \frac{\upi}{\uai} + \frac{2 \upi}{\ga} - \frac{2 \upi}{\uai} \nonumber \\
   &= \paren{\frac{\ga}{\uai} + \frac{\uai}{\ga / 2} - 2} \pdi.
\label{eq:pwd1}
\end{align}
At the later time, the same calculation gives
\begin{equation}
\pwdii = \paren{\frac{\ga}{\uaii} + \frac{\uaii}{\ga / 2} - 2} \pdii.
\label{eq:pwd2}
\end{equation}

\begin{figure}
\begin{centre}
\includegraphics[width=15cm]{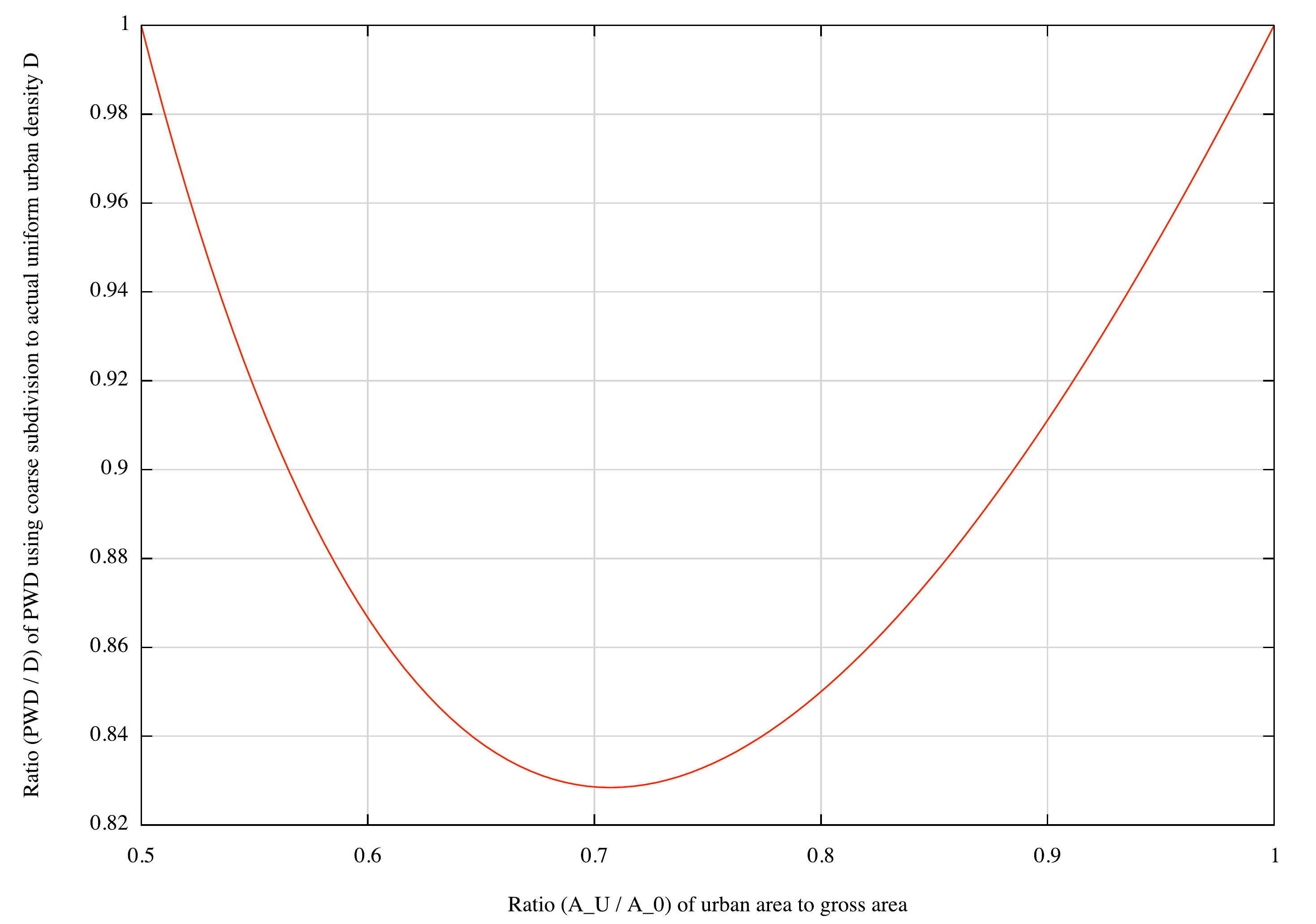}
\end{centre}
\caption{Variation of $\pwd / \pd$ as urban area \ua\ increases from $\ga / 2$ to \ga.}
\label{fig:coarse}
\end{figure}

It is of interest to observe the behaviour of the quantity in parentheses in both (\ref{eq:pwd1}) and (\ref{eq:pwd2}) as \ua\ increases from $\ga / 2$ to \ga.
This quantity is plotted in Figure \ref{fig:coarse}.
In the case $\ua = \ga / 2$, the inner parcel contains the entire population and the outer parcel no-one at all; the quantity in parentheses equals 1 and the PWD is the same as the uniform urban density \pd.
In the case $\ua = \ga$, on the other hand, both the inner and outer parcels contain equal populations at density \pd; the quantity in parentheses is again 1 and $\pwd = \pd$ again.

On the other hand, as \ua\ increases continuously from $\ga / 2$ to \ga, the quantity in parentheses actually \emph{declines} to a minimum at $\ga / \sqrt{2} \approx 0.7 \ga$, where it equals not 1 but $2 (\sqrt{2} - 1) \approx 0.828$.
It then increases back to the value 1 at $\ua = \ga$.

The implications of this for measurements of PWD over time in the above example are striking.
A case in point is where \uai\ is very slightly greater than $\ga / 2$, so that the urban area initially extends only a small way into the outer parcel, and \pwdi\ is only slightly less than the actual urban density \pdi\ at this time.
If the city subsequently expands so that \uaii\ is about 70\% of \ga, but the overall density of the urbanised area remains the same (so that $\pdi = \pdii$), then PWD may decline by as much as 17 per cent due to the expansion alone, \emph{even though} there has been no change in the actual density experienced by city residents.
Even more striking, if the density of the urbanised area increases over this time, say by 10 per cent (so that $\pdii = 1.1 \pdi$), then PWD could \emph{still} decline by almost 9 per cent, despite the increased density!

It should be clear that the example above could have been constructed with more than two parcels.
In particular, each of the two parcels could be subdivided into $N$ smaller parcels, and exactly the same results would be obtained, provided only that each of the $N$ `outer' parcels includes a similar proportion of urbanised land.
If the urbanised area \ua\ is roughly circular, for example, then one may subdivide the `outer' parcels in a wedge-shaped fashion (as indicated by the radial grid lines in Figure \ref{fig:schematic}), without affecting the end result.
Note that this type of subdivision with an outer `ring' of land parcels at similar distance from the city centre is typical of `urban interface' local government areas in cities like Melbourne.

This example illustrates a hazard with basing PWD calculations on coarse subdivisions of a region that includes both urban and non-urban land.
It may in particular help explain some counterintuitive trends in calculations of population density for Melbourne, summarised in Table \ref{tab:melbpd} based on figures from the Australian Bureau of Statistics (ABS).
\begin{table}[t]
\begin{centre}
\begin{tabular}{c|ccc|c}
& \multicolumn{2}{c}{Melbourne Urban Centre} & Overall density & \pwd\ using \\
Year & \up & \ua\ (ha) & \pd\ (p/ha) & SA3 parcels \\ \hline
1976 & 2,479,225 & 148,000 & 16.75 & \\
1981 & 2,578,759 & 162,000 & 15.90 & 17.2 \\
1986 & 2,645,484 & 165,000 & 16.05 & 16.0 \\
1991 & 2,761,995 & 166,500 & 16.60 & 15.4 \\
1996 & 2,932,149 & 179,100 & 16.40 & 15.2 \\
2001 & 3,160,171 & 208,030 & 15.20 & 15.3 \\
2006 & 3,371,889 & 215,280 & 15.65 & 15.8 \\
2011 & 3,707,530 & 254,320 & 14.60 & 16.7
\end{tabular}
\end{centre}
\caption{Overall density calculations for Melbourne Urban Centre,\protect\footnotemark\
compared with PWD for Melbourne Statistical Division on ABS `SA3' geography}
\label{tab:melbpd}
\end{table}
\footnotetext{%
Overall density figures up to 1991 sourced from Manning \cite{manning1984} and Mees \cite{mees2000}.
Subsequent OD figures and all PWD figures sourced by the author from ABS Census data.}
Focussing on the period 1981--91, it has been observed that PWD in Melbourne based on the relatively coarse SA3 geography (roughly equivalent to post-amalgamation local government areas) declined by about 10\% at this time; yet calculations by Manning \cite{manning1984} and Mees \cite{mees2000} of overall density based on the Melbourne `Urban Centre' (which draws a notional boundary around the actual urbanised area) show a modest increase of 5--6\% over the same period.
This was a time of consolidation within already-established middle and fringe suburbs coupled with a decline in population and hence density in inner suburbs, which resembles in its effects the simpler urban growth example given.

Figures for the period after 1991 illustrate a different problem with calculating overall density based on the ABS declared Urban Centre.
On several occasions the area \ua\ is seen to increase more rapidly than the urban population \up; this reflects assimilation to the urban area of large tracts of formerly extra-urban settlements bridged by substantial areas of rural land.
(ABS criteria allow for rural land parcels to be reclassified as urban if they separate urbanised regions up to 3km apart.)
This effect was particularly strong between 1976 and 1981, between 1996 and 2001 and between 2006 and 2011; the effect is to give a potentially misleading suggestion of declining urban density when in fact the trend since at least 1991 has been toward increasing consolidation.

To further illustrate the underlying trends as they play out through both the OD and PWD measures, Table \ref{tab:melbpdfix} shows the same indices for \emph{fixed} regions within the established Melbourne urban area.
\begin{table}
\begin{centre}
\begin{tabular}{c|ccc|ccc}
& \multicolumn{3}{c|}{Inner Melbourne} & \multicolumn{3}{c}{1981 established area} \\
& \multicolumn{3}{c|}{($\ua = 45060$ha)} & \multicolumn{3}{c}{($\ua = 66400$ha)} \\
Year & \up & \pd\ (p/ha) & PWD & \up & \pd\ (p/ha) & PWD \\ \hline
1981 & 1,097,937 & 24.4 & 27.5 & 1,558,215 & 23.5 & 25.9 \\
1986 & 1,068,102 & 23.7 & 26.4 & 1,511,632 & 22.8 & 24.9 \\
1991 & 1,064,725 & 23.6 & 26.1 & 1,501,031 & 22.6 & 24.5 \\
1996 & 1,081,273 & 24.0 & 26.5 & 1,504,762 & 22.7 & 24.8 \\
2001 & 1,118,911 & 24.8 & 27.0 & 1,546,100 & 23.3 & 25.2 \\
2006 & 1,190,476 & 26.4 & 28.7 & 1,627,257 & 24.5 & 26.6 \\
2011 & 1,283,802 & 28.5 & 31.0 & 1,752,995 & 26.4 & 28.7
\end{tabular}
\end{centre}
\caption{Overall density and PWD calculations for fixed established areas in Melbourne.\protect\footnotemark\
PWD calculations use ABS `SA3' geography.}
\label{tab:melbpdfix}
\end{table}
The consolidating trend is now evident in both the OD and PWD numbers, and unlike those for the entire urban area, the two measures move almost perfectly in lockstep.
As with the figures in Table \ref{tab:melbpd} however, there is no great difference between the PWD and OD measures and it may appear there is little to distinguish the two in practice.
But this is in turn largely a consequence of calculating PWD based on land parcels that mask the contrast in density at the neighbourhood level---a phenomenon discussed in more detail in the next section.

\section{The Perils of Large Parcels II: Sensitivity to Boundaries}

The concept of PWD is aimed at helping to overcome the inevitable `lumpiness' of population distribution in urban regions.
Some cities display significant clustering of population in high-density neighbourhoods or `urban villages' while others spread more uniformly over large areas, and it is useful to be able to compare the two from the point of view of `typical experienced density'.
PWD provides a useful step in this direction, but it still relies on subdividing the region at a sufficient level of detail to clearly delineate these clusters.

Where the division into parcels for the PWD calculation does not accurately distinguish higher-density clusters from lower-density neighbourhoods, the number that results can be unexpectedly sensitive to the often arbitrary placement of parcel boundaries, just as OD calculations are sensitive to the definition of the overall urban boundary.
\footnotetext{%
Figures sourced by the author from ABS Census data.
`Inner Melbourne' refers to the `Inner', `Inner East' and `Inner South' SA4 groupings as defined by ABS.
`1981 established area' refers to SA3 areas that were fully or almost fully urbanised in 1981.}

To demonstrate this, consider two parcels of land within a larger region, whose PWD is calculated from (\ref{eq:pwd}) as
\begin{equation}
\pwd_1 = \frac{1}{\gp} \paren{\frac{P_1^2}{A_1} + \frac{P_2^2}{A_2} + \text{other terms}}.
\label{eq:pwdsens1}
\end{equation}
Suppose that the boundary between parcels 1 and 2 is now perturbed very slightly, such that parcel 1 gains from parcel 2 an apartment block with $p$ residents, while parcel 2 gains from parcel 1 a pocket park covering the same (relatively small) area as the apartment block but with no residences.
The parcel areas $A_1$ and $A_2$ are therefore unchanged, but the populations are perturbed to $P_1 + p$ and $P_2 - p$ respectively.
The new PWD is calculated as
\begin{equation}
\begin{split}
\pwd_2 &= \frac{1}{\gp} \paren{\frac{(P_1 + p)^2}{A_1} + \frac{(P_2 - p)^2}{A_2} + \text{other terms}} \\
   &= \frac{1}{\gp} \paren{\frac{P_1^2}{A_1} + \frac{P_2^2}{A_2}
      + 2p \paren{\frac{P_1}{A_1} - \frac{P_2}{A_2}} + p^2 \paren{\frac{1}{A_1} + \frac{1}{A_2}}
      + \text{other terms}}.
\end{split}
\label{eq:pwdsens2}
\end{equation}
Comparing (\ref{eq:pwdsens2}) with (\ref{eq:pwdsens1}), it is seen that this minor change in boundary causes the overall PWD to change by an amount
\begin{equation}
\begin{split}
\pwd_2 - \pwd_1 &= 2 \frac{p}{\gp} \paren{\frac{P_1}{A_1} - \frac{P_2}{A_2}}
      + \frac{p^2}{\gp} \paren{\frac{1}{A_1} + \frac{1}{A_2}} \\
   &= 2 \frac{p}{\gp} \paren{\frac{P_1}{A_1} - \frac{P_2}{A_2} + \frac{p}{A_H}}
\end{split}
\label{eq:pwdsens}
\end{equation}
where $A_H$ is the harmonic mean of the parcel areas $A_1$ and $A_2$.
(If it happens that $A_1 = A_2$, then the quantity in parentheses is just $(P_1 - P_2 + p) / A_1$.)

What is the practical meaning of formula (\ref{eq:pwdsens})?
The first thing to note is that the boundary shift can only leave the PWD unchanged if parcel 2 started out slightly \emph{more} dense than parcel 1---specifically by the amount $p / A_H$.
Contrariwise, if the parcel densities $P_1 / A_1$ and $P_2 / A_2$ were initially the same, or parcel 1 was more dense than parcel 2, the PWD will always \emph{increase} as a result of the boundary shift.
This is despite there being no actual change `on the ground'.

It is not difficult to see from this that when one considers different alternative ways of dividing the \emph{same} urban region into the \emph{same} number of parcels, a much higher PWD number can result when at least some of the parcel boundaries are tightly drawn around the perimeter of higher-density residential clusters, than if the parcel boundaries are drawn simply to delineate roughly equal areas or to follow local administrative boundaries.

\begin{figure}
\begin{centre}
\includegraphics[width=7cm]{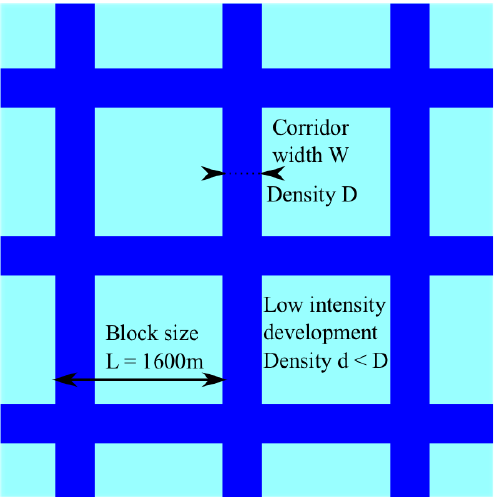}
\end{centre}
\caption{Stylised example of urban area with selective corridor development.}
\label{fig:cville}
\end{figure}

The effect can be significant even when the land parcels are on a `neighbourhood' scale aimed at getting a `fine-grained' analysis.
An example demonstrating this is depicted in Figure \ref{fig:cville}.
This hypothetical urban area has high-density development, at some uniform population density $D$, confined to the immediate vicinity of an arterial road grid spaced at intervals $L$, approximately 1.6km (one mile) apart.
Traditional suburban detached housing, at a lower uniform density $d$, has been maintained in the interior of the blocks delineated by the arterial grid.
The actual width of the high-density corridors is denoted $W$, and assumed to be less than $L / 2$.
While this development pattern has been idealised for purposes of discussion, it is reminiscent of the pattern emerging in established inner suburbs of cities like Melbourne.

\begin{figure}
\begin{centre}
(a)\includegraphics{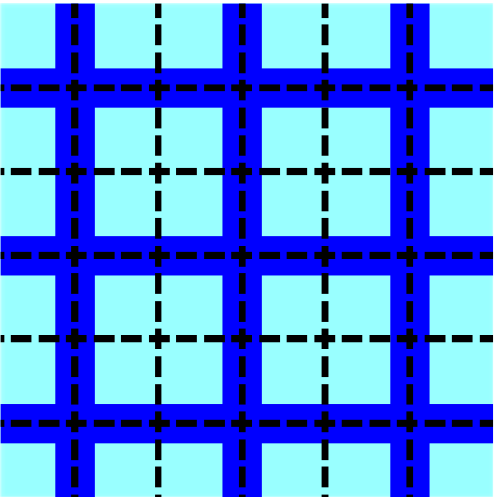}\hspace{2em}
(b)\includegraphics{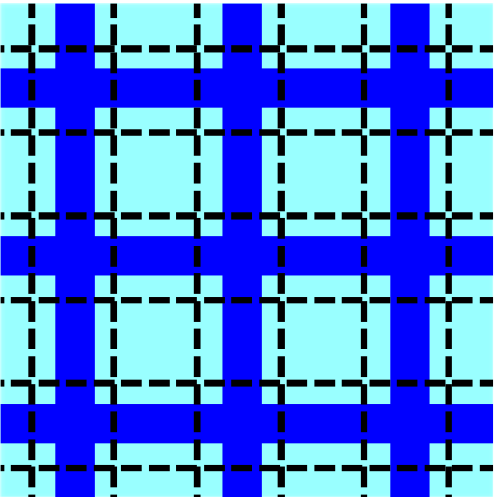}
\end{centre}
\caption{Two possible ways to subdivide into equal size parcels.}
\label{fig:subdivs}
\end{figure}

Now suppose the PWD of this neighbourhood is assessed by subdividing into square parcels of side length $L/2 = 0.8$km (half a mile).
Figure \ref{fig:subdivs}(a) shows one obvious way this subdivision might be carried out, with parcel boundaries aligned with the arterial roads.
The alternative subdivision in Figure \ref{fig:subdivs}(b) is very similar: the parcels are the same size and shape, but offset so that they are  \emph{centred} on the arterial roads---and accordingly on the high-density corridors.

These two choices of subdivision have radically different consequences for the PWD calculation.
Observe first that when the parcels are defined as in Figure \ref{fig:subdivs}(a), each one encloses a vertical and horizontal strip of width $W/2$ at density $D$ with the remaining area at density $d$.
The population of each parcel is therefore%
\footnote{Assuming $D$ and $d$ are expressed as persons per hectare and $L$ and $W$ in kilometres, a scaling factor of 100 is required to convert square kilometres to hectares.
A neater though more esoteric way to carry out these calculations in SI units is to express $L$ and $W$ in \emph{hectometres}: 1hm = 100 metres, 1km = 10hm and $1\mathrm{hm}^2 = 1$ hectare.}
\begin{equation}
\begin{split}
P_a &= 100 D \paren{2 \paren{\frac{W}{2}} \paren{\frac{L}{2} - \frac{W}{2}} + \paren{\frac{W}{2}}^2}
      + 100 d \paren{\frac{L}{2} - \frac{W}{2}}^2 \\
   &= 25 d L^2 + 25 \paren{D - d} W \paren{2 L - W}.
\end{split}
\label{eq:popsuba}
\end{equation}
Since all parcels in Figure \ref{fig:subdivs}(a) have an identical population distribution, the PWD of the entire region is equal to the average density of a single parcel, namely $P_a / A$ where $A = 100 (L/2)^2$ is the parcel area in hectares:
\begin{equation}
\pwd_a = \frac{P_a}{100 (L / 2)^2} = \frac{P_a}{25 L^2}
   = d + (D - d) \paren{\frac{W}{L}} \paren{2 - \frac{W}{L}}.
\label{eq:pwdsuba}
\end{equation}
Note that this quantity $\pwd_a$ is in fact the same as the overall density \od\ of the region, as a consequence of Proposition \ref{prop:pwdpd} when all parcels have the same density.

For the subdivision in Figure \ref{fig:subdivs}(b) the calculation is more involved, because the population distributions are no longer identical across all parcels.
Given the assumption $W \leq L / 2$, three broad types may be distinguished:
\begin{itemise}
\item
`Interior block' parcels with no high-density development, hence with uniform density $d$ and population
\begin{equation}
P_1 = 100 d \paren{\frac{L}{2}}^2 = 25 d L^2.
\label{eq:popsubb1}
\end{equation}
\item
`Arterial crossroad' parcels traversed by high-density development in both north-south and east-west directions, and population
\begin{equation}
\begin{split}
P_2 &= 100 D \paren{W^2 + 2 W \paren{\frac{L}{2} - W}}
      + 100 d \paren{4 \paren{\frac{L}{4} - \frac{W}{2}}^2} \\
   &= 25 d L^2 + 100 \paren{D - d} W \paren{L - W}.
\end{split}
\label{eq:popsubb2}
\end{equation}
\item
`Arterial transverse' parcels with one central high-density strip, either north-south or east-west, and population
\begin{equation}
P_3 = 100 D W \paren{\frac{L}{2}} + 100 d \paren{\frac{L}{2} - W} \paren{\frac{L}{2}}
   = 25 d L^2 + 50 \paren{D - d} W L.
\label{eq:popsubb3}
\end{equation}
\end{itemise}
The PWD of the region is that of a representative $2 \times 2$ square of parcels, comprising one parcel of the first type, one of the second and two of the third.
Given the area of each parcel is $A = 100 (L / 2)^2$ hectares, one calculates
\begin{equation}
\begin{split}
\pwd_b &= \frac{1}{P_1 + P_2 + 2 P_3} \paren{P_1 \frac{P_1}{A} + P_2 \frac{P_2}{A}
      + 2 P_3 \frac{P_3}{A}} \\
   &= \frac{P_1^2 + P_2^2 + 2 P_3^2}{100 \paren{d L^2 + \paren{D - d} W \paren{2 L - W}}
      \cdot 100 \paren{L / 2}^2} \\
   &= d \cdot \frac{1 + 2 \paren{D / d - 1} \paren{W / L} \paren{2 - W / L}
         + 2 \paren{D / d - 1}^2 \paren{W / L}^2 \paren{1 + 2 (1 - W / L)^2}}
      {1 + \paren{D / d - 1} \paren{W / L} \paren{2 - W / L}}.
\end{split}
\label{eq:pwdsubb}
\end{equation}

To get a better insight into the formulae (\ref{eq:pwdsuba}) and (\ref{eq:pwdsubb}) for the PWD according to the two different subdivisions, it helps to notice that they essentially depend on just two dimensionless factors:
\begin{itemise}
\item
the ratio $D / d$ of intense development density to `ambient' neighbourhood density; and
\item
the geometric ratio $W / L$ of corridor width to block size.
\end{itemise}
In addition, the PWD formulae themselves suggest defining a \emph{consolidation factor}
\begin{equation}
K = \paren{\frac{D}{d} - 1} \paren{\frac{W}{L}} \paren{2 - \frac{W}{L}}.
\label{eq:kconsol}
\end{equation}
In terms of this factor $K$, the PWD according to subdivision \ref{fig:subdivs}(a) is simply
\begin{equation}
\pwd_a = \od = d \paren{1 + K}.
\label{eq:pwdsubak}
\end{equation}
The factor $K$ therefore denotes the relative increase in overall density OD brought about by developing the high-density corridors.

Now, consider how formula (\ref{eq:pwdsubb}) might be simplified through the use of $K$.
The denominator of this formula works out to be just $1 + K$, while the numerator is $1 + 2K + Q$, where $Q$ is a formula resembling $K^2$.
Analysing the numerator as $(1 + K)^2 + (Q - K^2)$, and cancelling a factor $(1 + K)$, there results the formula
\begin{equation}
\pwd_b = d \paren{1 + K
   + \frac{\paren{D / d - 1}^2 \paren{W / L}^2 \paren{2 (1 - W / L)^2 + (W / L)^2}}{1 + K}}.
\label{eq:pwdsubbk}
\end{equation}
The calculated densities $\pwd_a$ and $\pwd_b$ may now be compared directly: their difference is
\begin{equation}
\pwd_b - \pwd_a = \frac{d}{1 + K} \paren{\frac{D}{d} - 1}^2 \paren{\frac{W}{L}}^2 \paren{2 \paren{1 - \frac{W}{L}}^2 + \paren{\frac{W}{L}}^2}.
\label{eq:pwdsubdiff}
\end{equation}
Notice that provided $D > d$ and $W > 0$ one always has $\pwd_b > \pwd_a = \od$, in accordance with Proposition \ref{prop:pwdpd}.
But more importantly, given realistic values of $D / d$ and $W / L$ the difference in (\ref{eq:pwdsubdiff}) is not only positive but substantial.
Consider for example the following realistic scenario:
\begin{itemise}
\item
$L = 1.6$km, as above (one-mile blocks are traditional in many cities);
\item
$W = 0.2$km (intense development extends 100 metres on either side of the main roads);
\item
$d = 15$ persons per hectare (typical for low density Melbourne suburbs); and
\item
$D = 195$ persons per hectare (in the typical range for 3--4 storey apartments).
\end{itemise}
This scenario leads to $K = 2.81$ (to three significant figures).
The PWD according to Figure \ref{fig:subdivs}(a), which is also the overall density, is $\pwd_a = \od = 57.2$ persons per hectare---an almost fourfold increase on the base density $d$.
The difference $\pwd_b - \pwd_a$, however, is given by formula (\ref{eq:pwdsubdiff}) as 13.7, so that the PWD according to Figure \ref{fig:subdivs}(b) is $\pwd_b = 70.9$ persons per hectare---nearly 25 per cent greater than $\pwd_a$ and a near fivefold increase on $d$.

\begin{figure}[t]
\begin{centre}
\includegraphics[width=\textwidth]{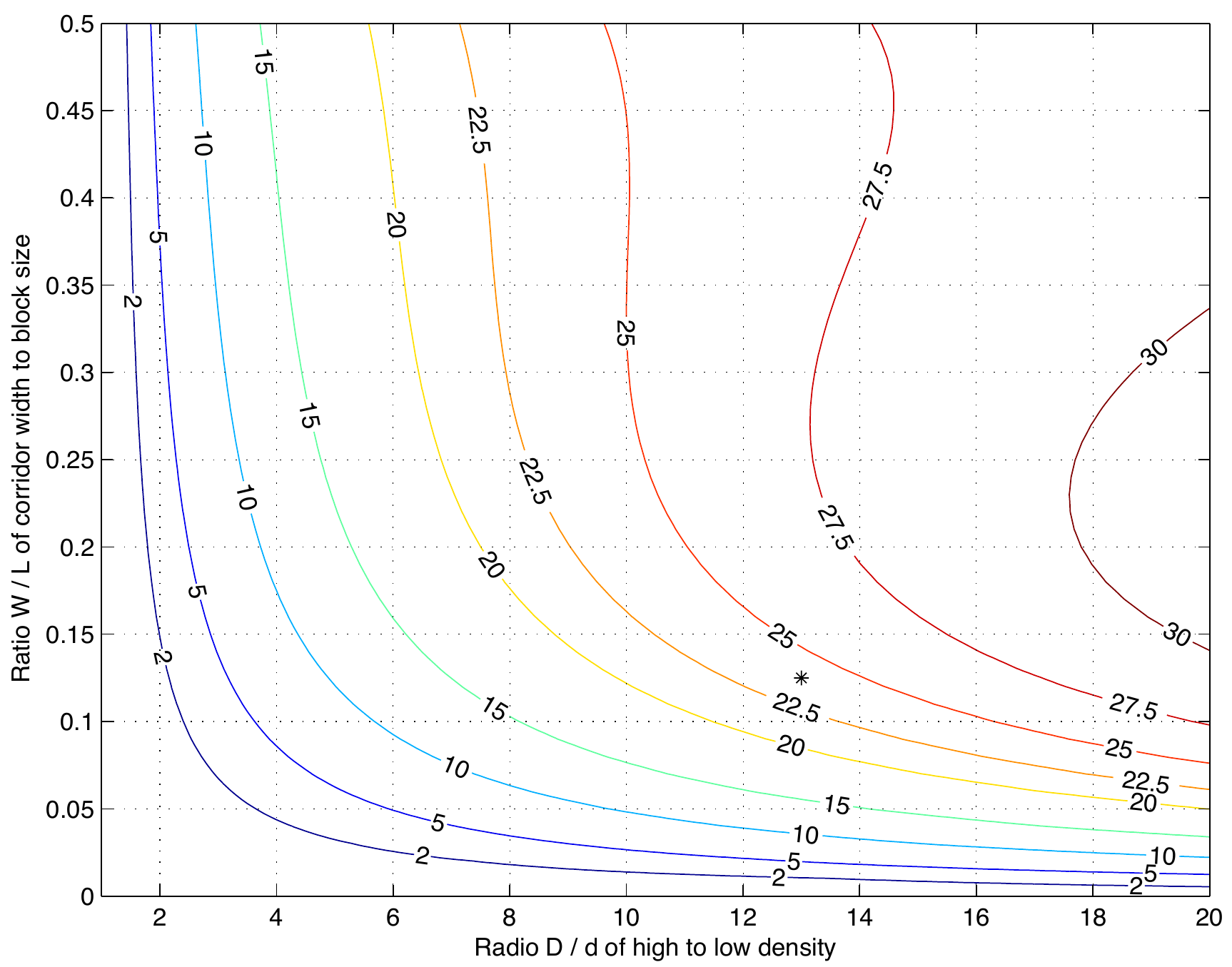}
\end{centre}
\caption{Contour plot of percentage difference $(\pwd_b - \pwd_a) / \pwd_a \times 100\%$ in PWD calculation for the two subdivisions in Figure \ref{fig:subdivs}.
An asterisk marks the scenario in the main text.}
\label{fig:pwdsubdiff}
\end{figure}

In general, Figure \ref{fig:pwdsubdiff} charts the percentage difference between the PWD calculated for the subdivisions in Figure \ref{fig:subdivs}(a) and \ref{fig:subdivs}(b), as a function of the two dimensionless factors $D / d$ and $W / L$.
The scenario above with $D / d = 13$ and $W / L = 0.125$ is marked with an asterisk on this chart.

Even larger variations in PWD are possible when the parcel boundaries are free to shift in other ways.
\begin{figure}
\begin{centre}
\includegraphics{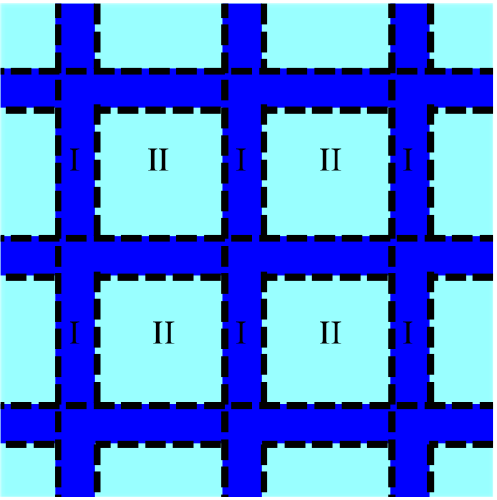}
\end{centre}
\caption{A subdivision that maximises the PWD for the region in Figure \ref{fig:cville}.}
\label{fig:subdivc}
\end{figure}
The largest PWD figure obtainable from the idealised urban area of Figure \ref{fig:cville} occurs when all the high-density area belongs to one parcel (or set of parcels) and all the low-density area to another.
A subdivision of this type is shown in Figure \ref{fig:subdivc}, where each block of size $L \times L$ is divided into a high-density L-shaped parcel (labelled I) and a low-density square parcel (labelled II).

For the subdivision in Figure \ref{fig:subdivc}, the population of a parcel of type I is
\begin{equation}
P_I = 100 D \paren{2 W (L - W) + W^2} = 100 D W (2L - W)
\label{eq:popsubc1}
\end{equation}
and that of a parcel of type II is
\begin{equation}
P_{II} = 100 d (L - W)^2.
\label{eq:popsubc2}
\end{equation}
The PWD of the region subdivided as in Figure \ref{fig:subdivc} is equivalent to that of a representative $L \times L$ block containing one parcel of each type---with density $D$ and $d$ respectively---and accordingly may be calculated as
\begin{equation}
\begin{split}
\pwd_c &= \frac{P_I D + P_{II} d}{P_I + P_{II}}
   = d \cdot \frac{\paren{D / d}^2 \paren{W / L} \paren{2 - W / L} + \paren{1 - W / L}^2}
      {\paren{D / d} \paren{W / L} \paren{2 - W / L} + \paren{1 - W / L}^2} \\
   &= d \cdot \frac{1 + \paren{D / d + 1} K}{1 + K}
   = d + \frac{K}{1 + K} D
\end{split}
\label{eq:pwdsubc}
\end{equation}
where $K$ is given by formula (\ref{eq:kconsol}).

For the specific scenario given above the effect of a subdivision as in Figure \ref{fig:subdivc} is remarkable.
With $d = 15$ and $D = 195$ persons per hectare, and $K$ equal to 2.81, formula (\ref{eq:pwdsubc}) gives a PWD of 159 persons per hectare, compared with 57.2 and 70.9 for the two subdivisions in Figure \ref{fig:subdivs}.
Not only is $\pwd_c$ more than 2.7 times greater than the overall density; it is also more than double the value $\pwd_b$ calculated from a subdivision into $L / 2 \times L / 2$ parcels approximately aligned with the high-density corridors.
The subdivision in Figure \ref{fig:subdivc} contains fewer parcels than those in Figure \ref{fig:subdivs}, but by drawing boundaries tightly around areas of high density yields a much higher PWD figure.

\begin{figure}
\begin{centre}
\includegraphics[width=\textwidth]{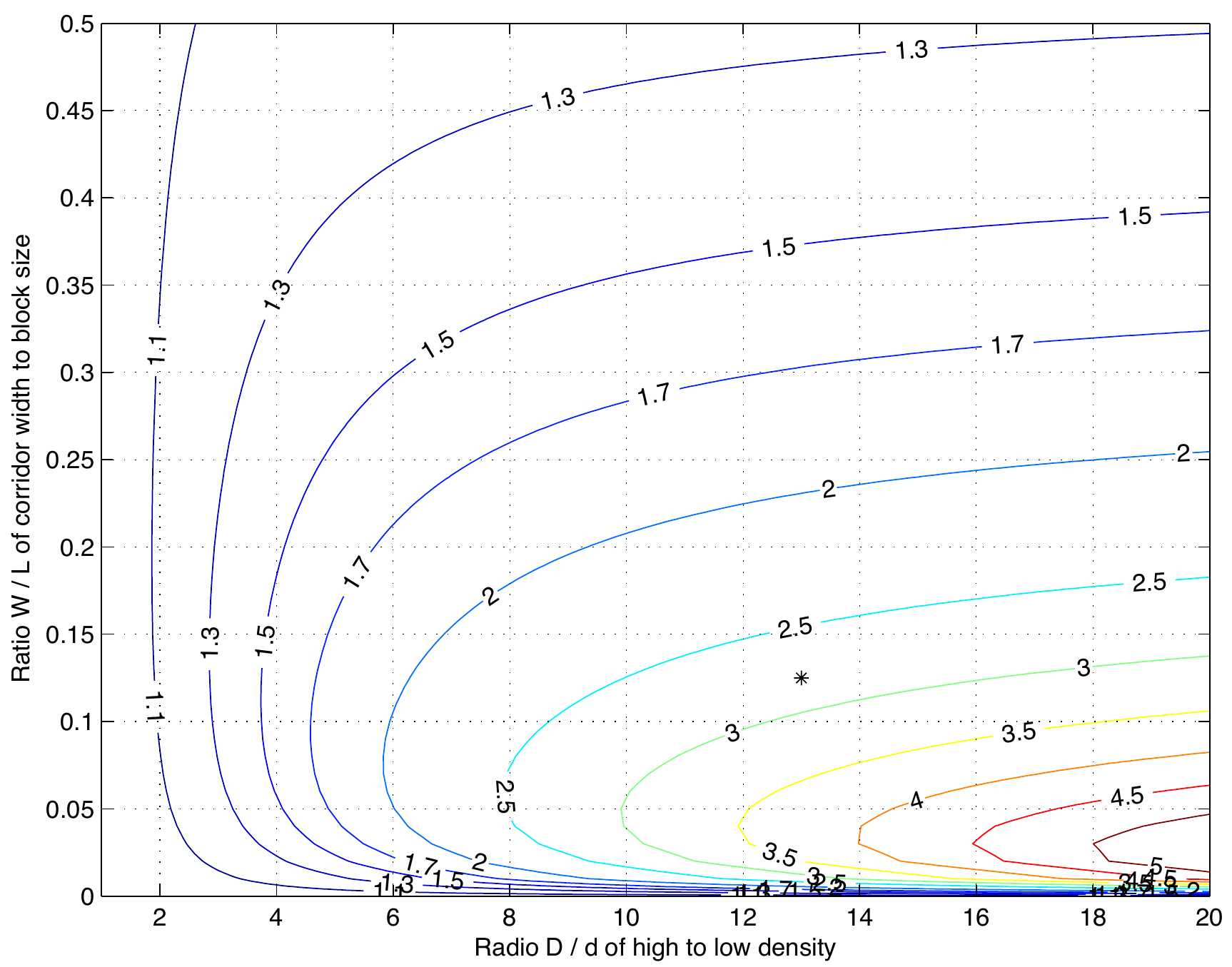}
\end{centre}
\caption{Contour plot of the ratio $\pwd_c / \pwd_a$ ($= \pwd_c / \od$) for the subdivision in Figure \ref{fig:subdivc}.
An asterisk marks the scenario in the main text.}
\label{fig:pwdsubcmul}
\end{figure}

Figure \ref{fig:pwdsubcmul} charts the ratio $\pwd_c / \pwd_a = \pwd_c / \od$ for the subdivision in Figure \ref{fig:subdivc} in the general case, as a function of the factors $D / d$ and $W / L$.
Again, the scenario given above is marked with an asterisk on this chart.
A notable feature here is that the increase in PWD is most pronounced when the high-density corridors are relatively narrow---and particularly so when the development density $D$ increases to compensate.

These results show that for a typical developing urban area with a patchwork of high and low density development, calculated values of PWD even with relatively small land parcels can vary over a wide range, and are sensitive to even small details of parcel boundaries.
Depending on how parcels are delineated the PWD number can be as low as the overall density or equal to several times that number.
In this context, obtaining a reliable figure for PWD is likely to invove adopting one of two alternative strategies:
\begin{enumerate}
\item
Ensure that parcels are small enough that any high-density development is contained entirely within a single parcel, rather than part of a larger parcel encompassing both high and low density areas.
(Note that the calculation leading to (\ref{eq:pwdsubc}) is unaltered if the parcels I and II are arbitrarily subdivided further, as long as each of the smaller parcels clearly belongs either to a type I or a type II region.)
\item
Alternatively, keep the parcel sizes moderate but design the parcel boundaries carefully to ensure that all areas with housing density significantly higher than in surrounding areas are tightly bounded.
\end{enumerate}
Either strategy, applied to a real-world neighbourhood with features reminiscent of Figure \ref{fig:cville}, ensures that the PWD obtained will be as close as practically possible to a value like $\pwd_c$---distinguished as being the limiting value obtained from arbitrarily fine-grained subdivisions of the study area.

\section{Conclusion}

Population-weighted density provides useful information about the population distribution in an urban area, beyond that conveyed by overall density or other traditional measures.
However, it is by no means immune to unexpected or seemingly paradoxical features, which can lead to errors of interpretation if one is not careful.

The PWD measure will always lead to a greater density number than the overall density, unless the region in question is absolutely uniform in density.
Indeed, the greater the extent to which population `clusters' within an urban area, the greater will be the amount by which PWD exceeds the overall density.
(This is what Bradford calls `clumpiness', and Eidlin the `density gradient index').
It has been shown that while PWD is sometimes referred to as ``giving equal weight to persons rather than hectares'', it is perhaps better characterised as a measure of clustering, whose effect is to give people an \emph{unequal} weighting based on the relative density of their neighbourhood.
Thus, PWD is closely related mathematically to the notion of a density-weighted population.

Particularly when the aim is to assess longitudinal trends in density within a fixed urban area, it has been found advantageous to subdivide the region as finely as is practical and to draw tight boundaries around areas of higher-density development.
When the subdivision is either too coarse or too heedless of development patterns, so that some of the subdivided parcels contain significant amounts of both urban and non-urban land, or both high-density and low-density development, there is a likelihood of paradoxical results.
In particular, the PWD of a growing region can appear to decline over a period even when the urbanised part of that region has increased in density; and the PWD even of a small mixed-density neighbourhood can be highly sensitive to the geometry of the land parcels used as input data.

In the end it is as Ernest Fooks might have said: to get a reliable, quantitative picture of density and `sprawl', one has to be prepared to X-ray one's city.
The thought experiments in this note, and the original Los Angeles paradox itself, should underline the fact that urban development can proceed in many more ways than can be neatly summed up with terms like `high density' or `sprawl'.
Misconceptions can arise from clinging to developmental stereotypes---such as the idea of ever-decreasing density as one proceeds from the core to the outskirts of cities---and it may be that there is no single number that neatly characterises `sprawl'.

Planners should instead be open to the idea that what at first looks like `sprawl' may reflect not the land use pattern as such, but rather the interaction of land use, the provision of urban transport or other public goods, and `urban geometry' factors independent of density.
In these situations one must look to the supporting infrastructure and services as the means to improved urban amenity.

\bibliographystyle{plain}
\bibliography{transport}

\end{document}